\documentclass[12pt]{amsart}
\usepackage[scale=0.81]{geometry} 		
\usepackage{graphicx}				
\usepackage[utf8]{inputenc}
\usepackage[english]{babel}
\usepackage[utf8]{inputenc}
\usepackage{float}
\usepackage{amsmath,amsfonts,amssymb,verbatim}
\usepackage{color}
\newcommand{\diag}{{\rm diag}}
\newcommand{\E}{{\rm E}}
\newcommand{\pen}{{\rm pen}}
\newcommand{\JKL}{{\rm JKL}}
\newcommand{\KL}{{\rm KL}}

\usepackage{amsthm}

\newtheorem{defi}{Definition}[section]
\newtheorem{theo}[defi]{Theorem}

\newtheorem{rmk}{Remark}

\newtheorem{assum}[defi]{Assumption}
\newcommand{\emilieFeb}[1]{\textcolor{black}{#1}}

\newcommand{\CG}[1]{\textcolor{black}{#1}}
\usepackage{xspace}

\providecommand{\keywords}[1]{\textbf{\textit{Keywords---}} #1}
\DeclareMathAlphabet\mathbfcal{OMS}{cmsy}{b}{n}
\usepackage{hyperref}
\hypersetup{
 citecolor=black,
 filecolor=black,
 linkcolor=black,
 urlcolor=black
}

\usepackage[T1]{fontenc}
\usepackage[utf8]{inputenc}
\usepackage{float}
\usepackage{setspace}
\usepackage{amsfonts,amsthm,mathrsfs}
\usepackage{psfrag,epsf}
\usepackage{enumerate}
\usepackage{url} 

\usepackage[english]{babel}
\usepackage{amsmath}

\begin{document}
\title{Joint rank and variable selection for parsimonious estimation in a \CG{high-dimensional} finite mixture regression model}

\author{Emilie Devijver}\address{Department of Mathematics and Leuven Statistics Research Center (LStat), KU Leuven, Leuven, Belgium}

\maketitle
\begin{abstract}
We study a dimensionality reduction \CG{technique} for finite \CG{mixtures} of high-dimensional multivariate response regression models. \CG{Both the dimension} of the response and the number of predictors \CG{are allowed to} exceed the
sample size. We consider predictor selection and rank reduction to obtain \CG{lower-dimensional} approximations.
A class of estimators with a fast rate of convergence is introduced.
We apply this result \CG{to} a specific procedure, introduced in \cite{procedures}, where the relevant predictors are selected by the Group-Lasso. 
\end{abstract}
\keywords{\CG{Finite mixture regression  ; Low rank estimation ; Non-asymptotic penalized criterion ; Oracle inequality ; Variable selection}


\tableofcontents

 \section{Introduction}
 The popular multivariate response regression model $\mathbf{Y}=\boldsymbol{\beta} \mathbf{X}+\boldsymbol{\varepsilon}$
postulates a linear relationship between \CG{a $q \times n$ matrix} $\mathbf{Y}$ containing $q$ responses for $n$ subjects and \CG{a $p \times n$ matrix $\mathbf{X}$ on $p$ predictor variables. 
The responses for individual $i \in \{ 1, \ldots, n\}$ are stored in the vector $\mathbf Y_i \in \mathbb{R}^q$ and the corresponding predictors are stored in the vector $\mathbf X_i \in \mathbb{R}^p$. In this model, the} error term $\boldsymbol{\varepsilon}$ is an unobserved $q \times n$ matrix with independent columns \CG{$\boldsymbol \varepsilon_1, \ldots , \boldsymbol \varepsilon_n\sim \mathcal{N}_q(0,\boldsymbol{\Sigma})$}.
The $q \times p$ regression matrix $\boldsymbol{\beta}$ and the $q \times q$ covariance matrix $\boldsymbol{\Sigma}$ \CG{are parameters that} need to be estimated. \CG{In practice, however,} this model is often too simple \CG{for} heterogeneous data. 

\CG{This paper considers} a finite mixture of $K$ linear multivariate regression models.
If the subject $i$ originates from the class  $k \in \{1,\ldots,K\}$, \CG{we assume that} there exists an unknown matrix of coefficients $\boldsymbol{\beta}_k\in \mathbb{R}^{q \times p}$ 
and an unknown covariance matrix  $\boldsymbol{\Sigma}_k \in \mathbb{S}_q^{++}$ such that the target variable \CG{$\mathbf{Y}_i$} results from a regression model
$
\mathbf{Y}_i=\boldsymbol{\beta}_k \mathbf{X}_i + \boldsymbol{\varepsilon}_i$, with $\boldsymbol\varepsilon_i \sim \mathcal{N}_q(0,\boldsymbol{\Sigma}_k)$.
If we use this model to study complex data with a high number of responses and predictors, 
the number of parameters could be quickly  larger than the sample size.
In order to deal with this issue, we propose to work with parsimonious models combining two well-known approaches.

The first approach consists in selecting relevant variables.
The Lasso estimator, introduced by Tibshirani \cite{Tibshirani}, is a classical tool.
In multiple multivariate regression models, it seems reasonable to take into \CG{account} the matrix structure. 
For \CG{this} purpose, we focus on row sparsity using the Group-Lasso estimator, introduced by Yuan and Lin \cite{YuanLin}. In our context, the groups correspond to the rows. 
\CG{Thus if, for $k \in \{ 1, \ldots, K\}$}, a matrix $\boldsymbol{\beta}_k$ has $|J|$ relevant rows, there are $|J|q$ coefficients to estimate instead of $pq$ per cluster. If $|J| \ll  p$,  the number of parameters to estimate \CG{is then drastically reduced}. 
Moreover, by selecting only the few relevant rows to describe the response, row sparsity may enhance interpretation.

The second approach is based on rank sparse models, introduced by Anderson \cite{Anderson}.
A linear model is rank sparse if the regression matrix has low rank or can at least  be well approximated by a low-rank matrix.
In our model, for \CG{each $k \in \{ 1, \ldots, K\}$}, if the matrix $\boldsymbol{\beta}_k$ has rank $\mathbf{R}_k$, $\boldsymbol{\beta}_k$ is fully determined by $\mathbf{R}_k(p+q-\mathbf{R}_k)$ coefficients,
which may be  smaller than the sample size $nq$.
These low-rank regression matrices appear in many applications\CG{, e.g.,} analysis of fMRI image data \cite{Harrison} \CG{and the} analysis of EEG data decoding \cite{Anderson}.
Moreover, \CG{note} that low-rank estimation generalizes principal component analysis, \CG{often used to} reduce the dimension of multivariate data.

Bunea et al. \cite{Bunea} mix these two estimators to provide a row and rank sparse estimator  for the linear model.
In practice,  row and rank sparse approaches are combined to benefit from the two methods: \CG{fewer} parameters to estimate, better interpretation thanks to the selection of relevant variables, and better estimation for data with low-rank structure.
In \cite{procedures}, a method has been proposed, called \CG{the} Lasso-Rank procedure, which combines both sparsity constraints for finite-mixture regression models. Heterogeneous data are then described by several clusters.

In this article, we focus on the theoretical aspects of combining row and rank sparse structures. \CG{Assuming such a structure} for linear models, Bunea et al. \cite{Bunea} obtain substantial oracle rate improvements over the existing methods, especially when the estimated rank is low. In this paper, we give a similar result for finite-mixture regression models. \CG{Note that, in contrast to Bunea et al.} \cite{Bunea}, we \CG{do not need to assume that} the covariance matrices \CG{$\boldsymbol{\Sigma}_1, \ldots, \boldsymbol{\Sigma}_K$ are} known. \CG{Therefore,} our result generalizes \CG{theirs} even \CG{in the case} $K=1$. Furthermore, our result \CG{does} not assume strong assumptions on the design and on the true \CG{underlying} density. \CG{We discuss its merits under classical assumptions and show that it can also serve to justify the model selection step of the Lasso-Rank procedure used to construct a collection of models.}

For the general model, we construct a penalty satisfying a non asymptotic control of the risk of the model minimizing the penalized likelihood, among a random collection of row and rank sparse models.
However, the penalty achieving the oracle type inequality is \CG{only} known up to multiplicative constants.
\CG{As these} quantities depend on bounds on \CG{the} parameters \CG{that} are difficult to construct \CG{in practice, we propose to calibrate the required constants} using slope \CG{heuristics}.

Since the geometry of mixture models is complex, we use tools which are different from \cite{Bunea}.
Our result \CG{exploits} a generalization of a  model selection theorem for maximum likelihood \CG{estimation} in a regression framework developed by Cohen and Le Pennec in \cite{CohenLePennec}.
The main difficulty \CG{lies} in the randomness of the collection of models \CG{under consideration}.
Indeed, as the number of responses and predictors is potentially high, the whole collection of models is not tractable in practice.
Thus, we \CG{do not look} at every possible sparsity level but \CG{rather} focus on a random subset constructed from the data.
As an example, in the Lasso-Rank procedure, this random subset corresponding to the relevant rows is constructed using the Group-Lasso estimator along the regularization path.
To solve this issue, we use an extension of the model selection theorem proved in \cite{inegOracleLassoMLE}.

The paper is organized as follows. 
In Section~\ref{sec:2}, we introduce finite \CG{mixtures} of linear regression models, and we describe the collection of models we consider. In Section~\ref{sec:3}, we define the divergences we need, and we derive an oracle type inequality.
Section~\ref{sec:4} investigates the Lasso-Rank procedure, which is a specific procedure used to construct the random subcollection. All proofs are \CG{relegated to an Appendix}.

\section{The model and the collection of models}
\label{sec:2}

\subsection{Mixture of linear regression models}
\label{mixtureModels}

\CG{Let} $(\mathbf{X}_1,\mathbf{Y}_1),\ldots,(\mathbf{X}_n,\mathbf{Y}_n)$ \CG{be a random sample, and let} 
 $\mathbf{x}_i$ and $\mathbf{y}_i$ \CG{denote the observed values of the} random variables $\mathbf{X}_i \in \mathbb{R}^p$ and $\mathbf{Y}_i \in  \mathbb{R}^q$, \CG{respectively}. \CG{We assume that} the response variable $\mathbf{Y}_i$ depends on the set of explanatory variables $\mathbf{X}_i$ through a regression-type model; \CG{its} conditional probability density function, denoted $s^\star$, 
is approximated by a mixture of $K$ multivariate Gaussian regression models.

If $(\mathbf{X}_i,\mathbf{Y}_i)$ originates from an individual $i$ in class $k$, we assume that
 there exists \CG{a matrix of class-specific regression coefficients} $\boldsymbol{\beta}_k \in \mathbb{R}^{q \times p}$,  
\CG{a $q\times q$ positive-definite matrix} $\boldsymbol{\Sigma}_k$, and a random vector $\boldsymbol \varepsilon_i \sim \mathcal{N}_q(0,\boldsymbol{\Sigma}_k)$, such that 
 $$
 \mathbf{Y}_i= \boldsymbol{\beta}_k \mathbf{X}_i +\boldsymbol \varepsilon_i.
 $$
We further assume that \CG{whatever the values of $\mathbf{x}_1, \ldots, \mathbf{x}_n$, the variables $\mathbf{Y}_1 | \mathbf{X}_1 = \mathbf{x}_1, \ldots,
\mathbf{Y}_n | \mathbf{X}_n = \mathbf{x}_n$ are mutually independent. For each $i \in \{ 1, \ldots , n\}$, we take} $\mathbf{Y}_i|\mathbf{X}_i=\mathbf{x}_i \sim s_{\xi}^K(\mathbf{y}|\mathbf{x}_i)$, with
  \begin{align}
  \label{modele}
  &s_{\xi}^K(\mathbf{y}|\mathbf{x})=\sum_{k=1}^{K} \frac{\boldsymbol{\pi}_k}{(2 \pi)^{q/2}  \text{det}(\boldsymbol{\Sigma}_k)^{1/2}} \exp \left\{ -\frac{(\mathbf{y}-\boldsymbol{\beta}_{k} \mathbf{x}){^\top} \boldsymbol{\Sigma}_{k}^{-1}(\mathbf{y}-\boldsymbol{\beta}_{k} \mathbf{x})}{2} \right\},
   \end{align}
where $\xi=(\boldsymbol{\pi}_1,\ldots, \boldsymbol{\pi}_K,\boldsymbol{\beta}_{1},\ldots,\boldsymbol{\beta}_{K},\boldsymbol{\Sigma}_{1},\ldots,\boldsymbol{\Sigma}_{K}) \in \Xi_K$ with
\begin{equation}
\label{XiK}
\Xi_K = \Pi_{K} \times (\mathbb{R}^{q\times p})^K\times (\mathbb{S}_q^{++})^K.
\end{equation}

\CG{In the above, $\boldsymbol{\pi}_k$ denotes the proportion for} class $k \in \{1,\ldots,K\}$. For each such class and all $z \in \{1,\ldots,q\}$ and $\mathbf{x} = (\zeta_1,\ldots,\zeta_p) \in \mathbb{R}^p$, $[\boldsymbol{\beta}_{k} \mathbf{x}]_z= \sum_{j=1}^{p} [\boldsymbol{\beta}_k]_{z,j} \zeta_{j}$ is the $z$th component of the $k$th mixture component. \CG{We set} $ \Pi_{K} = \{ (\boldsymbol{\pi}_1, \ldots, \boldsymbol{\pi}_K) : \boldsymbol{\pi}_1 + \cdots + \boldsymbol{\pi}_K = 1 \} \subset (0,1)^K$ and let $
 \mathbb{S}_q^{++}$ \CG{denote the collection} of symmetric positive definite matrices on $\mathbb{R}^q$.

Moreover, we assume that, for all $k \in \{1,\ldots, K\}$, the covariance matrices $\boldsymbol{\Sigma}_k$ are diagonal. 
As $q$ is potentially large, it reduces drastically the number of parameters to estimate.
\CG{From} a theoretical viewpoint, if $q$ is not too large, the current result works for full covariance matrices, focusing on the eigenvalues. However, if $q$ is large, we are \CG{unable} to estimate these covariance matrices \CG{unless we make some additional} assumptions.
One way to allow for correlations might be to construct block-diagonal covariance matrices as in \cite{DevijverGallopin}. \CG{This strategy, however, will not be further considered here.}

\subsection{Framework and collection of models} 

In \CG{what follows}, a row will be said to be {relevant} if every coefficient of the row is equal to $0$ in every cluster.
\CG{Otherwise, the row is said to be irrelevant}. 
We denote by $J$ the set of indices of relevant rows.

\CG{Given a} model with $K$ clusters, to combine rank and row sparsity, we consider models with rank $\mathbf{R} = (\mathbf{R}_1,\ldots,\mathbf{R}_K)$ and relevant rows indexed by $j \in J$.
We denote by ${S}_{(K,J,\mathbf{R})}$ the set of conditional probability density functions corresponding to a mixture of $K$ regression models, when focusing on $|J|$ relevant predictors indexed by $j \in J$, and with matrix of regression
 coefficients $\boldsymbol{\beta}_k$ of rank $\mathbf{R}_k$ in class $k$, viz.
\begin{align}
 {S}_{(K,J,\mathbf{R})} = \bigl\{ s_{\xi}^{(K,J,\mathbf{R})} \bigr\},
\end{align}
where 
\begin{multline*}
s_{\xi}^{(K,J,\mathbf{R})} (\mathbf{y}|\mathbf{x}) =\sum_{k=1}^K \frac{\boldsymbol{\pi}_k }{(2 \pi)^{q/2} \det (\boldsymbol{\Sigma}_k)^{1/2}} \\ \times
\exp \left[ -\frac{1}{2}\{\mathbf{y}-(\boldsymbol{\beta}_k^{\mathbf{R}_k})^{[J]} \mathbf{x}\}{^\top} \boldsymbol{\Sigma}_k^{-1} \{\mathbf{y}-(\boldsymbol{\beta}_k^{\mathbf{R}_k})^{[J]}\mathbf{x}\} \right]
\end{multline*}
with
$\xi = (\boldsymbol{\pi}_1,\ldots,\boldsymbol{\pi}_K, \boldsymbol{\beta}_1^{\mathbf{R}_1},\ldots, \boldsymbol{\beta}_K^{\mathbf{R}_K},\boldsymbol{\Sigma}_1,\ldots,\boldsymbol{\Sigma}_K) \in \Xi_{(K,J,\mathbf{R})}$, $
\Xi_{(K,J,\mathbf{R})} = \Pi_K \times \Psi_{(K,J,\mathbf{R})}\times (\mathbb{S}_q^{++})^K$, and
\begin{multline*}
\Psi_{(K,J,\mathbf{R})} = \Bigl\{((\boldsymbol{\beta}_1^{\mathbf{R}_1})^{[J]},\ldots,(\boldsymbol{\beta}_K^{\mathbf{R}_K})^{[J]}) \in (\mathbb{R}^{q\times p})^K \text{ such that}\\
\forall_{k\in \{1,\ldots,K\}} \; \text{rank}(\boldsymbol{\beta}_k^{\mathbf{R}_k})=\mathbf{R}_k\Bigr\}.
\end{multline*}
\CG{Here,} $A^{[J]}$ \CG{stands for} the matrix $A$ \CG{whose} coefficients indexed by $j \in J$ are the nonzero coefficients. \CG{The dimension of the model $S_{(K,J,\mathbf{R})}$ is}
$$
D_{(K,J,\mathbf{R})} = \sum_{k=1}^K \{ \mathbf{R}_k (p+q-\mathbf{R}_k)+q+1\} -1.
$$

The complete collection of models \CG{consists of} models with $K$ clusters, row sparsity defined by $J$ describing all subsets of $\{1,\ldots,p\}$, 
and vector of ranks for a fixed $K$ 
and fixed $J$ belonging to $\{1,\ldots, |J|\wedge q\}^K$, viz.
\begin{align}
\label{collecModelS}
\mathcal{S} = \bigcup_{K \in \{1, 2, \ldots\}} \bigcup_{J \in \mathcal{P}(\{1,\ldots,p\})} \bigcup_{\mathbf{R} \in \{1,\ldots,|J|\wedge q\}^K} {S}_{(K,J,R)} = \bigcup_{(K,J, \mathbf{R}) \in {\mathcal{M}}} {S}_{(K,J,R)}.
\end{align}

This collection is generally large and hence not tractable in practice.
\CG{Therefore, we} bound the number of clusters by \CG{restricting our attention to a finite set $\mathcal{K} \subset \{ 1, 2, \ldots \} $}. We further focus on a \CG{(potentially random)} subcollection $\tilde{\mathcal{J}}$ of $\mathcal{P}(\{1,\ldots,p\})$ 
\CG{of} moderate size, and we reduce the number of possible vectors
 of ranks considered by introducing $\mathcal{R}_{(K,J)} \subset \{1\ldots,|J|\wedge q\}^K$, \CG{where in general, $a \wedge b = \min (a,b)$.
Accordingly, we consider}
 \begin{align}
 \label{modeleStilde}
\tilde{\mathcal{S}} = \bigcup_{K \in \mathcal{K}} \bigcup_{J \in \tilde{\mathcal{J}}} \bigcup_{\mathbf{R} \in \mathcal{R}_{(K,J)}} {S}_{(K,J,R)} = \bigcup_{(K,J, \mathbf{R}) \in \tilde{\mathcal{M}}} {S}_{(K,J,R)}.  
 \end{align}

\CG{In order to derive theoretical results, we further} assume that \CG{the} parameters are bounded. \CG{Our bounds are expressed in terms of} the singular values of $\boldsymbol{\beta}_k$ and \CG{the} eigenvalues of the covariance matrices \CG{$\boldsymbol{\Sigma}_1, \ldots, \boldsymbol{\Sigma}_K$}.
If we fix $K$, the number of components, the relevant columns $J$ and the vector of ranks $\mathbf{R}$, we define the bounded model by
\begin{multline}
\label{ensBorneSec4}
 {S}_{(K,J,\mathbf{R})}^\mathcal{B}=\bigg\{ s_{\xi}^{(K,J,\mathbf{R})} \in {S}_{(K,J,\mathbf{R})} :  \xi = (\boldsymbol{\pi},{\boldsymbol{\beta}},{\boldsymbol{\Sigma}}), \text{ such that}  \\
\forall_{k \in \{1,\ldots,K\}} \; \boldsymbol{\Sigma}_k = \diag([\boldsymbol{\Sigma}_{k}]_{1,1},\ldots,[\boldsymbol{\Sigma}_{k}]_{q,q}), 
\forall_{z \in \{1,\ldots,q\}} \; a_{\boldsymbol{\Sigma} }\leq [\boldsymbol{\Sigma}_{k}]_{z,z} \leq A_{\boldsymbol{\Sigma}},\\
\forall_{k \in \{1,\ldots,K\}} \; \boldsymbol{\beta}_k^{\mathbf{R}_k} =\sum_{r=1}^{\mathbf{R}_k} [\sigma_{k}]_r [u_k{^\top}]_{\bullet,r} [v_k]_{r,\bullet}, 
\forall_{r \in \{1,\ldots,\mathbf{R}_k\}} \; [\sigma_{k}]_r < A_\sigma \bigg\}.
\end{multline}
\CG{In the above,} for $k \in \{1,\ldots,K \}$, $[\sigma_k]_1, \ldots, [\sigma_k]_{\mathbf{R}_k}$ denote the singular values of $\boldsymbol{\beta}_k^{\mathbf{R}_k}$, with corresponding unit vectors $u_k$ and $v_k$. \CG{We are then led} to consider the following collection of models:
\begin{align}
\label{collectionModeles}
\tilde{\mathcal{S}}^{\mathcal{B}} = \bigcup_{K \in \mathcal{K}} \bigcup_{J \in \tilde{ \mathcal{J}}  }\bigcup_{\mathbf{R} \in \mathcal{R}_{(K,J)}} {S}_{(K,J,\mathbf{R})}^\mathcal{B}. 
\end{align}
\CG{Finally,} we assume that \CG{the} covariates $\mathbf{X}$ belong to an hypercube. \CG{Without loss of generality, this hypercube is taken to be} $ \in [0,1]^p$.

\section{Oracle inequality}
\label{sec:3}

In this section, we first introduce \CG{the Kullback--Leibler} divergence and some extensions \CG{thereof that are useful} to compare two conditional probability density functions. \CG{We then} derive an oracle type inequality which ensures that if we have penalized the log-likelihood in \CG{an appropriate} way, 
we  select a model which is as good as the oracle.

\label{oracle}

\subsection{{\CG{Kullback--Leibler} and \CG{Jensen--Kullback--Leibler} divergences}}
\label{divergences}

We first need to define divergences used to compare conditional probability density functions. The \CG{Kullback--Leibler} divergence is defined, \CG{for any two probability density functions $s$ and $t$,} by
\begin{equation}
{\KL(s,t)} = \begin{cases} \int_{\mathbb{R}^q} \ln \left\{ \frac{s(y)}{t(y)}\right\} s(y) dy & \mbox{if } sdy \ll  tdy,\\
\infty & \mbox{otherwise.} \end{cases}
\end{equation}
To deal with regression data, for observed covariates $\mathbf{x}_1,\ldots, \mathbf{x}_n$, we define
\begin{equation}
\KL^{\otimes_n} (s,t)=\E\left[\frac{1}{n} \sum_{i=1}^n \KL\{s( \cdot |\mathbf{x}_i),t( \cdot |\mathbf{x}_i)\} \right]
\end{equation}
for \CG{any two conditional probability density functions $s$ and $t$.}

\CG{For any two probability density functions $s$ and $t$,} we also define the \CG{Jensen--Kullback--Leibler} divergence, first introduced in \cite{CohenLePennec}, by
$$
\JKL_{\rho} (s, t)=\frac{1}{\rho} \, \KL\{s,(1-\rho)s+\rho t\}
$$
for any $\rho \in (0,1)$. The tensorized \CG{version} is defined by 
$$
\JKL_{\rho}^{\otimes_n} (s,t) = \E\left[\frac{1}{n} \sum_{i=1}^n \JKL_{\rho}\left\{s( \cdot |\mathbf{x}_i),t( \cdot |\mathbf{x}_i)\right\}\right]
$$
for \CG{any two conditional probability density functions $s$ and $t$.}
We use the \CG{Jensen--Kullback--Leibler} divergence rather than the \CG{Kullback--Leibler} divergence because it is bounded.
This boundedness turns out to be crucial to control the loss of the penalized maximum likelihood estimator under mild assumptions on the complexity of the model and on the parameters.

\CG{Note that the above} divergences \CG{are not metrics} because they do not satisfy the triangular inequality and they are not symmetric. \CG{Nevertheless,} they are \CG{very useful for comparing} two probability density functions.

\subsection{Oracle inequality for a general procedure to select relevant columns} 

Among the collection of models $\tilde{\mathcal{S}}^{\mathcal{B}}$ defined in \eqref{collectionModeles}, with various ranks and various level of sparsities, 
we want to select a model which is close to the best one.
The oracle is by definition the model belonging to the collection which minimizes the contrast with the true model (here, the contrast considered is the log-likelihood).

In practice, we \CG{cannot identify} the true model and hence the oracle is unknown. Nevertheless, we could select a model close to the oracle. More precisely, we aim at selecting, among $\tilde{\mathcal{S}}^\mathcal{B}$, the optimal set of predictors $J$, the optimal number of clusters $K$, and the optimal vector of ranks $\mathbf{R}$.

First, for each model defined by $(K,J,\mathbf{R}) \in \tilde{\mathcal{M}}$, we consider the conditional probability density function $\hat{s}^{(K,J,\mathbf{R})}$ where the parameters are estimated by the maximum likelihood estimator.
Among all $(K,J,\mathbf{R}) \in \tilde{\mathcal{M}}$, we want to select the density $\hat{s}^{(\hat{K},\hat{ J}, \hat{\mathbf{R}})}$ which is the closest one to the true distribution $s^\star$ according to \CG{the} divergences defined in Section \ref{divergences}.

Even if the oracle is not \CG{attainable} in practice, \CG{it turns out that we can do almost as well as the oracle. This result, which is our main finding, is formally stated below and proved in the Appendix.}

\begin{theo}
\label{OracleInequality}
Let \CG{$(\mathbf{x}_1, \mathbf{y}_1), \ldots, (\mathbf{x}_n, \mathbf{y}_n) \in [0,1]^p \times \mathbb{R}^q$} 
be the observations arising from an unknown conditional probability density function $s^\star$.
 Let $\mathcal{M}$ denotes the set of every potential model and $\tilde{\mathcal{M}} $ the restricted set, as defined in \eqref{collecModelS} and \eqref{modeleStilde}, respectively.
 
\emilieFeb{ Assume \CG{that there exist} $\tau >0$ and $\delta_{\KL}>0$ such that, for all $(K,J, \mathbf{R}) \in \mathcal{M}$, \CG{one can find} $\bar{s}^{(K,J,\mathbf{R})} \in {S}_{(K,J,\mathbf{R})}^\mathcal{B}$ such that 
 \begin{equation}
 \label{hypTau}
 \KL^{\otimes_n} (s^\star,\bar{s}^{(K,J,\mathbf{R})}) \leq  \inf_{t \in {S}_{(K,J,\mathbf{R})}^{\mathcal{B}}} \KL^{\otimes_n}(s^\star,t) +\frac{\delta_{\KL}}{n} \hspace{0.5cm} \text{and} \hspace{0.5cm} \bar{s}^{(K,J,\mathbf{R})} \geq e^{-\tau} s^\star. 
 \end{equation}
 }
For $(K,J,\mathbf{R}) \in \mathcal{M}$, \CG{let} $\hat{s}^{(K,J,\mathbf{R})}$ \CG{denote} the maximum likelihood estimator for the model $S^\mathcal{B}_{(K,J,\mathbf{R})}$, satisfying
 \begin{align*}
  \hat{s}^{(K,J,\mathbf{R})} &= \underset{s_\xi^{(K,J,\mathbf{R})} \in {S}^\mathcal{B}_{(K,J,\mathbf{R})}}{\operatorname{argmin}} 
  \left[ -\frac{1}{n} \sum_{i=1}^n \ln\left\{ s_\xi^{(K,J,\mathbf{R})} (\mathbf{y}_i|\mathbf{x}_i) \right\} \right].
 \end{align*}
 Then \CG{the following statements hold true.}
 \begin{itemize}
 \item [a)] \CG{There exist} absolute constants $\kappa>0$ and $C$ such that whenever
  \begin{align}
&\pen(K,J,\mathbf{R}) \nonumber\\
&\geq \kappa  \frac{D_{(K,J,\mathbf{R})}}{n} \left[2 B^2(A_\sigma,A_{\boldsymbol{\Sigma}},a_{\boldsymbol{\Sigma}})
-\ln\left\{ \frac{D_{(K,J,\mathbf{R})}}{n} B^2(A_\sigma,A_{\boldsymbol{\Sigma}},a_{\boldsymbol{\Sigma}}) \wedge 1 \right\}\right.  \nonumber\\
 & \left.\hspace{2.5cm}+  (1 \vee \tau) \ln\left\{ \frac{4 epq}{(K(1+|J|)-q^2)\wedge pq} +\sum_{k=1}^K \mathbf{R}_k \right\}\right] \label{penalite}
 \end{align}
 for every $(K,J,\mathbf{R}) \in \mathcal{M}$, with  $B(A_\sigma,A_{\boldsymbol{\Sigma}},a_{\boldsymbol{\Sigma}})$  a constant depending only on parameter bounds,
  \CG{there exists a} random vector $(\hat{K},\hat{J},\hat{\mathbf{R}}) \in \tilde{\mathcal{M}}$ such that
\begin{align*}
 (\hat{K},\hat{J},\hat{\mathbf{R}}) &= \underset{(K,J,\mathbf{R}) \in \tilde{\mathcal{M}}}{\operatorname{argmin}} \left[-\frac{1}{n} \sum_{i=1}^n \ln \{\hat{s}^{(K,J,\mathbf{R})} (\mathbf{y}_i|\mathbf{x}_i)\}
 +\pen(K,J,\mathbf{R}) \right].
\end{align*}
\item [b)]
Whatever the true density $s^\star$,
\begin{multline}
\label{inegalite oracle proc}
 \E\bigl\{\JKL_{\rho}^ {\otimes_n} (s^\star,\hat{s}^{(\hat{K},\hat{J},\hat{\mathbf{R}})} )\bigr\} \\
 \leq C \E\Bigl[ \inf_{(K,J,\mathbf{R}) \in \tilde{\mathcal{M}}} \Bigl\{\inf_{t \in {S}_{(K,J,\mathbf{R})}^\mathcal{B}} \KL^ {\otimes_n } (s^\star,t) \\
 + \pen(K,J,\mathbf{R}) \Bigr\} + \frac{(1 \vee \tau) }{n} \Bigr].
\end{multline}
\end{itemize}
\end{theo}

\bigskip
This theorem is deduced from an adaptation of a model selection theorem for maximum likelihood estimators developed  by Massart  \cite{MassartStFlour}.
This adaptation, derived in \cite{procedures}, allows \CG{us} to focus on regression models (see also \cite{CohenLePennec}) and to focus on a random subcollection of the whole collection of models.

\CG{In Theorem 3.1, the critical assumptions} are the control of the bracketing entropy of each model in the whole collection of models and the construction of weights for each model to control 
the complexity of the collection of models. Assumption \eqref{hypTau} on the true density $s^\star$  is \CG{needed} because we consider a random subcollection of models $\tilde{\mathcal{S}}^\mathcal{B}$ from the whole collection of models $\mathcal{S}^\mathcal{B}$. 
Thanks to this assumption, we \CG{can invoke Bernstein's} inequality to control the additional randomness. The parameter $\tau$ depends on the true unknown density $s^\star$ and cannot be explicitly determined for this reason. Specific cases \CG{can be discussed}, e.g., the true conditional probability density function $s^\star$ should not be bounded along $\mathbf{Y}$ to allow \eqref{hypTau}. Moreover, if the true density is Gaussian, we need to bound the variance of the models \CG{from below}, which is related to the collection of models constructed in \eqref{ensBorneSec4}.

\CG{While we could make appropriate assumptions} on the true density $s^\star$ \CG{in order to determine $\tau$ explicitly,} we choose not to do so. Under the assumption that the \CG{Kullback--Leibler} divergence and the Hellinger distance are equivalent, we can explicitly determine $\tau$. 
Then, the result is still meaningful as the equivalence between \CG{the Kullback--Leibler} divergence and \CG{Hellinger's} distance is \CG{well known}.
Note that the parameter $\tau$ appears in the oracle type inequality \eqref{inegalite oracle proc}, but also in the penalty term \eqref{penalite}. 
We may construct a larger penalty independent of $\tau$, achieving an oracle-type inequality, but the rate will be larger. Details are discussed in \CG{the Appendix.}

For all $(K,J,\mathbf{R}) \in \mathcal{M}$, the mixture parameters are bounded in order to construct brackets over ${S}^\mathcal{B}_{(K,J,\mathbf{R})}$, and thus to bound the entropy \CG{from above}.

Because we consider a random sub-collection of models which is small enough, our estimator $\hat{s}^{(K,J,\mathbf{R})}$ is tractable in practice.
Moreover, we \CG{remark} that this is a \CG{finite sample} result, which allows us to study cases for which $p$ increases with $n$.

\section{Lasso-Rank procedure}
\label{sec:4}

In this section, we describe briefly the Lasso-Rank procedure introduced in \cite{procedures}, and we derive a corollary of Theorem \ref{OracleInequality} in this specific context.

\subsection{Description of the main steps}
The Lasso-Rank procedure is a  method used to construct a collection of models as defined in \eqref{modeleStilde}.
We recall here the main steps of this procedure. \CG{To this end, we assume that} we have access to a set of potential number of components $\mathcal{K}$.

\paragraph{A) Row selection}
Fix $K \in \mathcal{K}$. To detect the relevant rows, the empirical contrast is penalized by a group-penalty on the conditional mean parameters 
proportional to 
$$||\mathbf{P}_k^{-1} \boldsymbol{\beta}_k||_\text{gr} = \sum_{j=1}^p {\sqrt{q} \sqrt{ \sum_{z=1}^q [\mathbf{P}_k^{-1} \boldsymbol{\beta}_k]_{z,j}}},$$
where the Cholesky decomposition $\mathbf{P}_k{^\top} \mathbf{P}_k = \boldsymbol{\Sigma}_k^{-1}$ defines $\mathbf{P}_k$  for all $k \in \{1,\ldots,K\}$. 
We then  consider, for \CG{some choice of} regularization parameter $\lambda>0$, 
\begin{align}
\label{Lassosec4bis}
 \hat{\xi}^{\text{GL}} (\lambda) &= \underset{\xi=(\boldsymbol{\pi},\boldsymbol{\boldsymbol{\beta}},\boldsymbol{\boldsymbol{\Sigma}}) \in \Xi_K}{ \operatorname{argmin}}
 \left[ -\frac{1}{n} \sum_{i=1}^n \ln \{s_\xi^K(\mathbf{y}_i|\mathbf{x}_i)\}  + \lambda \sum_{k=1}^K \boldsymbol{\pi}_k ||\mathbf{P}_k^{-1} \boldsymbol{\beta}_k||_\text{gr}\right]
\end{align}
with $\Xi_K$ defined in \eqref{XiK}.

We remark that the penalty takes into account the mixture weight as, e.g., introduced in \cite{KhaliliChen}. 
This is in line with the common practice in relating the penalty to the sample size.
We refer the interested reader to \CG{Section 7.1.4 of} \cite{Stadler} for an interesting empirical discussion about a penalty similar to this one, focusing on this mixture weight.

To get models \CG{that are} more or less sparse, we compute the Group-Lasso estimator defined in  \eqref{Lassosec4bis} with several regularization parameters $\lambda$, 
and deduce \CG{sets of relevant rows}, denoted by $\hat{J}_{(K,\lambda)}$ for all $\lambda$.
We denote by $\mathcal{J}^{GL} = \cup_\lambda \hat{J}_{(K,\lambda)}$.

\begin{rmk}
\em
Notice that here, the Group-Lasso penalty is used to select a subset of variables for one choice of regularization parameter in the Lasso-Rank procedure. 
One could also use the Group-Lasso penalty to get a ranking of the variables, as done, e.g., in Bach \cite{Bach}.
\end{rmk}

\paragraph{B) Refitting}
 We estimate \CG{the} parameters by maximum likelihood estimation under rank constraint on the restricted  set of  relevant rows.

\paragraph{C) Model selection}

Varying $K \in \mathcal{K}$, $J \in \mathcal{J}^{\text{GL}}, \mathbf{R} \in \mathcal{R}_{(K,J)}$, we obtain a collection of models
$$\mathcal{S} = \bigcup_{ K \in \mathcal{K}} \bigcup_\mathcal{J \in {\mathcal{J}}^{\text{GL}} }\bigcup_{\mathbf{R} \in \mathcal{R}_{(K,J)}} {S}_{(K,J,\mathbf{R})}.$$
In practice, in \cite{procedures}, the slope heuristic is used, which is a data-driven non-asymptotic criterion.
It considers a penalized criterion, where the penalty is proportional to the dimension. There are two ways to calibrate this proportionality constant, but both consider the minimal coefficient
$\kappa_{\text{min}}$, the optimal one $\kappa_{\text{opt}}$ being approximated by twice the minimal coefficient $\kappa_{\text{min}}$.

One calibration method is the \textit{Slope Heuristic Dimension Jump} (SHDJ):
 $\kappa_\text{min}$ corresponds to the largest dimension jump on the graph representing the model dimension as a function of the coefficient $\kappa$. Another method is the
\textit{Slope Heuristic Robust Regression} (SHRR):  $\kappa_\text{min}$ corresponds to the slope of a 
robust regression performed between the log-likehood and the model dimension for complex models. 

The two methods are derived from the same heuristic and they offer two different visual checks of the adequacy of the model selection procedure to the data. They should select the same model. 
\CG{For an implementation of this approach, see Baudry et al.  \cite{baudry} and the \texttt{Capushe} package.}

\subsection{{Oracle inequality for relevant columns selected by the Group-Lasso estimator}}

In this section, we comment the use of Theorem \ref{OracleInequality} for this specific Lasso-Rank procedure. In general, we are not able to combine this theorem with results on the Group-Lasso to bound the oracle term, and then to bound the divergence between the model selected and the true conditional probability density function.
However, \CG{in some specific contexts}, we \CG{can} control the oracle.

For example, if $s^\star$ corresponds to a linear regression with sparse rank and \CG{a} few selected row, \CG{then} under some conditions on the design matrix, one can prove that this procedure is adaptive to row and rank sparsity.
This kind of result is discussed \CG{in Section~2 of Bunea et al.} \cite{Bunea}, and we refer to  Lounici et al. \cite{Lounici} for results on the Group-Lasso estimator.
However, our theorem stands non asymptotically and without \CG{any} assumptions on the design.
We have proved that we do the best possible \CG{within} the collection, \CG{although specific conditions are required} to prove that the best model \CG{within} the random collection is indeed \CG{objectively} good.
There is here a \CG{dirth} of theoretical tools \CG{and this challenging topic deserves further consideration.}

\subsection{Use of the slope heuristic}

\CG{From a theoretical point of view}, we want to ensure that the slope heuristic which penalizes the log-likelihood with the model dimension  selects a good model. We \CG{choose} a penalty 
$$
\text{pen}(K,J,R) = C \, \frac{D_{(K,J,R)}}{n} + \tilde{C} \, \frac{D_{(K,J,R)}}{n} \ln \left(\frac{D_{(K,J,R)}}{n}\right)
$$
with explicit constants $C$ and $\tilde{C}$ depending on absolute constants and on some bounds on the parameters. \CG{In practice, however, these bounds} are not tractable, and \CG{so} we prefer to calibrate $C$ and $\tilde{C}$ 
from the data using the slope heuristic. Moreover, we consider a simple version of the penalty term, viz.
$$
\text{pen}(K,J,R) = C \, \frac{D_{(K,J,R)}}{n}
$$
as already proposed \CG{in} \cite{DevijverGallopin} and \cite{Lebarbier}.
The calibration of $C$ from data has originally been proposed and \CG{validated} in the context of heteroscedastic regression with fixed design \cite{Baraud,BirgeMassart}, and for \CG{homoscedastic} regression with fixed design \cite{Arlot}.
In this paper, we proved that the shape of an optimal penalty is close to the one used in the slope heuristic, up to a \CG{logarithmic} term.

\bigskip
\noindent
\textbf{Acknowledgments.}
I am grateful to Pascal Massart and Christophe Giraud for stimulating discussions.
I also thank the referees, the Associate Editor and the \CG{Editor-in-Chief, Christian Genest,} for comments and suggestions which \CG{led to substantial improvements. This research was supported in part} by the Interuniversity Attraction Poles Program (IAP-network P7/06), Belgian Science Policy Office, \CG{whose funding} is gratefully acknowledged.

 \bigskip
 \noindent

\setcounter{section}{0}
\section*{Appendix}

\CG{This Appendix contains} the details of the proof of Theorem \ref{OracleInequality}. 
It is derived from a general model selection theorem stated in \CG{Part~A}.
The proof of Theorem \ref{OracleInequality} \CG{can be derived by checking} \CG{Assumptions} \ref{AsumptionHm},  \ref{AsumptionSepm} and \ref{AsumptionK} described in \CG{Parts~B and C}.

\bigskip
\noindent
\textbf{A. A general oracle type inequality for model selection}

\bigskip
Model selection \CG{can be performed using either the AIC or the BIC} criterion \cite{akaike,Schwarz}. \CG{An important limitation of these criteria, however,} is that they are asymptotic. \CG{Their use in high-dimensional settings is thus ad hoc}.
The slope heuristic is a method which overcomes this difficulty. It \CG{was} introduced by Birg\'e and Massart  \cite{BirgeMassart}.

In \cite{MassartStFlour}, Massart proposed a model selection theorem which guarantees that a penalized criterion leads to a good model selection, the penalty being defined by the model complexity.
It \CG{provides support} for the slope heuristic \CG{approach in a finite sample setting}.
Cohen and Le Pennec  \cite{CohenLePennec} generalized this theorem \CG{to a} regression framework.
\CG{Here, we need a} generalization of \CG{their result} detailed  in \cite{inegOracleLassoMLE} because we consider a random collection of models $(S_m)_{m \in \mathcal{M}}$, indexed by $\mathcal{M}$. 

Let $(\mathbf{X},\mathbf{Y}) \in \mathcal{{X}}\times \mathcal{{Y}}$.
We observe $\mathbf{x}_1, \ldots, \mathbf{x}_n$.
\CG{Let us start} with a description of the \CG{required} assumptions.
First, we impose a structural assumption.
It is a bracketing entropy condition on the model with respect to the tensorized Hellinger divergence 
$$
(d^{\otimes_n}_{H})^2(s,t) = \E \left[ \frac{1}{n} \sum_{i=1}^n d^2_{H} \{s(\cdot |\mathbf{x}_i),t(\cdot |\mathbf{x}_i)\} \right]
$$
for two conditional probability density functions $s$ and $t$.
A bracket $[\ell,u]$ is a pair of functions such that for all $(\mathbf x, \mathbf y) \in \mathcal{{X}} \times \mathcal{{Y}}$, 
$
\ell(\mathbf y,\mathbf x) \leq s(\mathbf y|\mathbf x) \leq u (\mathbf y,\mathbf x)$.
For $\varepsilon >0$, the bracketing entropy $\mathcal{H}_{[\cdot]} (\varepsilon,S,d_H^{\otimes_n})$ of a set $S$ is defined as the logarithm of the minimum number 
of brackets $[\ell,u]$ of width $d_H^{\otimes_n}(\ell,u)$ smaller than $\varepsilon$ such that every   probability density function of $S$ belongs to one of these brackets. \CG{In what follows, $m$ denotes any fixed element of $\mathcal{M}$.}

{\renewcommand{\thedefi}{$\text{H}_m$}
\begin{assum}
\label{AsumptionHm}
There exists a non-decreasing function $\psi_m$ such that 
$\varpi \mapsto 1/\varpi \psi_m (\varpi)$ is non-increasing on 
$(0, \infty)$. Furthermore, for all $\varpi \in \mathbb{R}^+$ and $s_m \in S_m$,
$$
\int _0^\varpi \sqrt{ \mathcal{H}_{[\cdot]} (\varepsilon,S_m(s_m,\varpi),d_H^{\otimes_n})} d \varepsilon \leq  \psi_{m}(\varpi) ,
$$
where $S_m(s_m,\varpi)=\{t \in S_m: d_H^{\otimes_n}(t,s_m) \leq \varpi \}$.
The model complexity $\mathcal{D}_m$ is then defined by $n \varpi^2_m$ with $\varpi_m$ \CG{standing for} the unique solution of 
\begin{align}
\label{model complexity}
\frac{1}{\varpi} \, \psi_{m}(\varpi) = \sqrt{n}\varpi.
\end{align}
\end{assum}
\addtocounter{defi}{-1}}

\CG{Observe} that the model complexity \CG{does not} depend on the bracketing entropies of the global models $S_m$, but \CG{rather on those of the} smaller localized sets $S_m(s_m,\varpi)$. \CG{This} is a weaker assumption.

For technical reasons, a separability assumption is also required. \CG{In what follows, $\lambda$ denotes Lebesgue's measure.}

{\renewcommand{\thedefi}{$\text{Sep}_m$}
\begin{assum}
\label{AsumptionSepm}
There exists a countable subset $S^{'}_m$ of $S_m$ and a set $\mathcal{{Y}}_m^{'}$ with $\lambda(\mathcal{{Y}} \setminus \mathcal{{Y}}^{'}_m)=0$ such that for every $t \in S_m$, there exists a sequence $(t_\ell)_{\ell \geq 1}$ of elements of $S_m^{'}$ such that for every $\mathbf x$ and every $\mathbf y \in \mathcal{{Y}}_m^{'}$,
\CG{$\ln\{t_\ell (\mathbf y|\mathbf x)\} \to \ln\{t(\mathbf y|\mathbf x)\}$ as $\ell \to \infty$}.
\end{assum}
\addtocounter{defi}{-1}}

\bigskip
\CG{This assumption allows us to} work with a countable subset. We also need an 
information theory type assumption on our collection of models. We assume the existence of a Kraft-type inequality for the 
collection.

{\renewcommand{\thedefi}{K}
\begin{assum}
\label{AsumptionK}
There exists a family $(w_m)_{m \in \mathcal{M}}$ of non-negative numbers such that
$$
\sum_{m \in \mathcal{M}} e^{-w_m} \leq \Omega <  \infty.
$$
 \end{assum}
\addtocounter{defi}{-1}}

\bigskip
\CG{We are now in a position to state our global theorem for a random subcollection of models in a regression framework.}

\bigskip
\noindent
\textbf{Theorem A}. {\em
 Assume we observe $(\mathbf{x}_1,\mathbf{y}_1), \ldots, (\mathbf{x}_n,\mathbf{y}_n) \in \mathcal{X}\times \mathcal{Y}$ with unknown conditional probability density function $s^\star$. Let $\mathcal{S} = (S_m)_{m\in \mathcal{M}}$
 be \CG{an} at most countable collection of conditional probability density function sets.
\CG{Suppose that} Assumption \ref{AsumptionK} holds and that Assumptions \ref{AsumptionHm} and \ref{AsumptionSepm} are valid for every \CG{model} $S_m \in \mathcal{S}$. \CG{Further assume that there} exist $\tau>0$ and $\delta_{\KL}>0$  such that for all $m \in \mathcal{M}$, \CG{there is an element} $\bar{s}_m$ of $S_m$ such that
 \begin{equation}
 \label{hyp tau}
 \emilieFeb{
 \KL^{\otimes_n} (s^\star,\bar{s}_m) \leq  \inf_{t \in S_m} \KL^{\otimes_n}(s^\star,t) +\frac{\delta_{\KL}}{n} \hspace{0.5cm} \text{and} \hspace{0.5cm} \bar{s}_m \geq e^{-\tau} s^\star.
 }
 \end{equation}

 Introduce $(S_m)_{m \in \tilde{\mathcal{M}}}$ some random sub-collection of $(S_m)_{m \in \mathcal{M}}$.
 For $m\in  \tilde{\mathcal{M}}$, consider the collection $\hat{s}_m$ of $\eta$-log-likelihood minimizer in $S_m$, satisfying
 $$\sum_{i=1}^n - \ln\{\hat{s}_m (\mathbf{y}_i|\mathbf{x}_i)\} \leq \inf_{s_m \in S_m} \left[ \sum_{i=1}^n - \ln\{s_m(\mathbf{y}_i|\mathbf{x}_i)\} \right] + \eta.$$
 
 Then, for any $\rho \in (0,1)$ and any $C_1 > 1$, there are two constants $\kappa_0$ and $C_2$ depending only on $\rho$ and $C_1$ such that, as soon as for every index $m \in \mathcal{M}$, 
 \begin{equation}
\label{penalite2}
 \pen(m) \geq \kappa \{\mathcal{D}_m + (1 \vee \tau) w_m\}
   \end{equation}
with $\kappa > \kappa_0$, and where the model complexity $\mathcal{D}_m$ is defined through the equation \eqref{model complexity},
 the penalized likelihood estimate $\hat{s}_{\hat{m}}$ with $\hat{m} \in \tilde{\mathcal{M}}$ such that
 $$\sum_{i=1}^n - \ln \{\hat{s}_{\hat{m}} (\mathbf{y}_i|\mathbf{x}_i)\} + \pen(\hat{m}) \leq \inf_{m \in \tilde{\mathcal{M}}} 
 \left[ \sum_{i=1}^n - \ln\{ \hat{s}_m (\mathbf{y}_i |\mathbf{x}_i)\} + \pen (m) \right] + \eta^{'}$$
 satisfies
\begin{align}
\label{inegalite oracle}
 \E\{\JKL_{\rho}^ {\otimes_n} (s^\star,\hat{s}_{\hat{m}} )\} &\leq 
 C_1 \E\left\{ \inf_{m \in \tilde{\mathcal{M}}} \inf_{t \in S_m} \KL^ {\otimes_n } (s^\star,t)+2 \frac{\pen(m)}{n} \right\} \\
 &+ C_2 (1 \vee \tau) \frac{\Omega^2}{n} + \frac{\eta' + \eta}{n}. \nonumber
\end{align}}

\CG{Under the conditions of this theorem, we are thus able to choose, within a random subcollection of models, a representative} which is as good as the oracle, up to a constant $C_1$, and some additive terms \CG{of order} $1/n$.
This result is non-asymptotic, and gives a theoretical penalty to select this model.

 The proof of this theorem is detailed in  \cite{inegOracleLassoMLE}. Nevertheless, we \CG{sketch} the main ideas \CG{here in order to give the reader a better grasp of} the assumptions.
 From Assumptions \ref{AsumptionHm} and \ref{AsumptionSepm}, \CG{one can call on} maximal inequalities to \CG{control}, except on a set \CG{whose} probability \CG{is} less than 
 $e^{-w_{m^{'}}-w}$ for all $w$,  the ratio of the centered empirical process of $\ln(\hat{s}_{m'})$ over the Hellinger distance between $s^\star$ and $\hat{s}_{m'}$; \CG{the bound is of the order of} $1/n$.
 Thanks to \CG{Bernstein's} inequality, \CG{which is valid in view of inequality \eqref{hyp tau}}, we \CG{thus} control the process over the random subcollection.
 Thanks to Assumption \ref{AsumptionK}, we  sum up each event over the collection of models.
 It leads to the oracle type inequality \eqref{inegalite oracle}.

Now, to prove Theorem \ref{OracleInequality}, we have to satisfy Assumptions \ref{AsumptionHm} and \ref{AsumptionK}. \CG{As for} Assumption \ref{AsumptionSepm}, \CG{it holds} for our conditional probability density functions.
 
 \bigskip
\noindent
\textbf{B. Assumption \ref{AsumptionHm}}

 \bigskip
\noindent
\textbf{B1. Decomposition of the mixture conditional probability density function}. \CG{Proceeding as in} \cite{CohenLePennec}, we can decompose the entropy by
\begin{align}
 \label{entropie}
\mathcal{H}_{[\cdot]} (\varepsilon,\mathcal{S}^{\mathcal{B}}_{(K,J,\mathbf{R})}, d_H^{\otimes_n}) &\leq \mathcal{H}_{[\cdot]} (\varepsilon,\Pi_K, d_H^{\otimes_n}) + 
\sum_{k=1}^K \mathcal{H}_{[\cdot]} (\varepsilon,\mathcal{F}_{(J,\mathbf{R}_k)}, d_H^{\otimes_n}),
\end{align}

where 
\begin{align*}
\mathcal{S}^{\mathcal{B}}_{(K,J,\mathbf{R})}&= \left\{
\begin{array}{lll}
&\mathbf y \in \mathbb{R}^q | \mathbf x \in \mathbb{R}^p \mapsto s_{\xi}^{(K,J,\mathbf{R})}(\mathbf y|\mathbf x), \\
&s_{\xi}^{(K,J,\mathbf{R})}(\mathbf y|\mathbf x)= \sum_{k=1}^K \boldsymbol{\pi}_k \varphi \left(\mathbf y | (\boldsymbol{\beta}_k^{\mathbf{R}_k})^{[J]} \mathbf x,\boldsymbol{\Sigma}_k\right), \\
&\xi = \left( \boldsymbol{\pi}_1,\ldots,\boldsymbol{\pi}_K, (\boldsymbol{\beta}_1^{\mathbf{R}_1})^{[J]},\ldots,(\boldsymbol{\beta}_K^{\mathbf{R}_K})^{[J]},\boldsymbol{\Sigma}_1,\ldots,\boldsymbol{\Sigma}_K \right),\\
&\xi \in \Xi_{(K,J,\mathbf{R})}  = \Pi_K \times \tilde{\Psi}_{(K,J,\mathbf{R})} \times ([a_{\boldsymbol{\Sigma}},A_{\boldsymbol{\Sigma}}]^q)^K\\
\end{array}
 \right\},
 \end{align*}
 $$
\Psi_{(K,J,\mathbf{R})} = \left\{((\boldsymbol{\beta}_1^{\mathbf{R}_1})^{[J]},\ldots,(\boldsymbol{\beta}_K^{\mathbf{R}_K})^{[J]}) \in (\mathbb{R}^{q\times p})^K : \text{rank}(\boldsymbol{\beta}_k)=\mathbf{R}_k\right\},
$$
\begin{multline*}
\tilde{\Psi}_{(K,J,\mathbf{R})} = \Bigg\lbrace ((\boldsymbol{\beta}_1^{\mathbf{R}_1})^{[J]},\ldots,(\boldsymbol{\beta}_K^{\mathbf{R}_K})^{[J]}) \in \Psi_{(K,J,\mathbf{R})} : \text{ for all } k\in \{1,\ldots,K\},\\
\boldsymbol{\beta}_k^{\mathbf{R}_k} = \sum_{r=1}^{\mathbf{R}_k} \sigma_r [\mathbf{u}]_{r,\bullet}{^\top}[\mathbf{v}]_{r,\bullet}, \text{ with } \sigma_r < A_\sigma  \text{ for all }r \in \{1,\ldots,\mathbf{R}_k\} \Bigg\rbrace,
\end{multline*}
$ \Pi_{K} = \{ (\boldsymbol{\pi}_1, \ldots, \boldsymbol{\pi}_K) : \boldsymbol{\pi}_1 + \cdots + \boldsymbol{\pi}_K = 1 \} \subset (0,1)^K$
and
\begin{multline*}
 \mathcal{F}_{(J,{R})}  =\left\{ \varphi(.|(\boldsymbol{\beta}^{R})^{[J]} \mathbf{X},\boldsymbol{\Sigma}) : \boldsymbol{\beta}^{R} = \sum_{r=1}^{R} \sigma_r [\mathbf{u}]_{r,\bullet}{^\top}[\mathbf{v}]_{r,\bullet}, \text{ with } \sigma_r < A_\sigma, \right. \\
  \left. \phantom{\sum^{1^{\frac{1}{2^5}}} }\boldsymbol{\Sigma} = \diag(\boldsymbol{\Sigma}_{1,1},\ldots,\boldsymbol{\Sigma}_{q,q})\in [a_{\boldsymbol{\Sigma}},A_{\boldsymbol{\Sigma}}]^q \right\},
\end{multline*}
where $ \varphi$ denotes the Gaussian density. For the proportions, it is known (see, e.g., \cite{Wasserman}) that
$$
\mathcal{H}_{[\cdot]} (\varepsilon,\Pi_K, d_H^{\otimes_n}) \leq \ln \left\{K (2 \pi e)^{K/2} \left( \frac{3}{\varepsilon}\right)^{K-1} \right\}.
$$

\bigskip
\noindent
\textbf{B.2 Entropy of the Gaussian   probability density functions with low rank.}
We need to bound the integrated entropy. \CG{To this end}, first we have to construct some brackets to recover $\mathcal{F}_{(J,R)}$ \CG{for fixed rank $R$ and arbitrary set $J$ of relevant rows}.
Fix $f \in \mathcal{F}_{(J,R)}$. We are looking for functions $\ell$ and $u$ such that $\ell \leq f \leq u$.
Because $f$ is Gaussian, we will choose $\ell$ and $u$ as \CG{Gaussian dilatations. Therefore, we need only} to determine the mean, the variance and the dilatation coefficients of $\ell$ and $u$.
\CG{To construct these brackets, we rely on the following two lemmas; see \cite{Maugis} for details}.

\bigskip
\noindent
\textbf{Lemma A}
{\em
 Let $\varphi( \cdot |\boldsymbol{\mu}_1,\boldsymbol{\Sigma}_1)$ and $\varphi( \cdot|\boldsymbol{\mu}_2,\boldsymbol{\Sigma}_2)$ be two Gaussian   probability density functions. \CG{Suppose that their covariance matrices are diagonal, with $ \boldsymbol{\Sigma}_a=\diag([\boldsymbol{\Sigma}_{a}]_{1,1},\ldots,[\boldsymbol{\Sigma}_a]_{q,q})$
 for $a \in \{1,2\}$, and that} $[\boldsymbol{\Sigma}_{2}]_{z,z}>[\boldsymbol{\Sigma}_1]_{z,z}>0$
 for all $z \in \{1,\ldots,q\}$. Then, for all $x \in \mathbb{R}^q$,
\begin{align*}
\frac{\varphi(x|\boldsymbol{\mu}_1,\boldsymbol{\Sigma}_1)}{\varphi(x|\boldsymbol{\mu}_2,\boldsymbol{\Sigma}_2)} 
\leq &\prod_{z=1}^q \frac{\sqrt{[\boldsymbol{\Sigma}_{2}]_{z,z}}}{\sqrt{[\boldsymbol{\Sigma}_{1}]_{z,z}}}\\
&\times e^{\frac{1}{2} (\boldsymbol{\mu}_1-\boldsymbol{\mu}_2){^\top} \diag\left( \frac{1}{[\boldsymbol{\Sigma}_{2}]_{1,1}-[\boldsymbol{\Sigma}_{1}]_{1,1}},\ldots,\frac{1}{[\boldsymbol{\Sigma}_{2}]_{q,q}-[\boldsymbol{\Sigma}_{1}]_{q,q}}\right) (\boldsymbol{\mu}_1-\boldsymbol{\mu}_2) }.
\end{align*}}

\bigskip
\noindent
\textbf{Lemma B}. {\em
The Hellinger distance \CG{between} two Gaussian probability density functions with diagonal covariance matrices is given by
\begin{multline*}
d_H^2\{\varphi( \cdot |\boldsymbol{\mu}_1,\boldsymbol{\Sigma}_1),\varphi( \cdot |\boldsymbol{\mu}_2,\boldsymbol{\Sigma}_2)\} 
 =  2-2 \left(\prod_{z=1}^{q}\frac{2 \sqrt{[\boldsymbol{\Sigma}_{1}]_{z, z}[\boldsymbol{\Sigma}_{2}]_{z, z}}}{[\boldsymbol{\Sigma}_{1}]_{z,z} + [\boldsymbol{\Sigma}_{2}]_{z,z}}\right)^{1/2} \\
 \times \exp\left[-\frac{1}{4} (\boldsymbol{\mu}_1-\boldsymbol{\mu}_2){^\top} \diag \left\{ \left(\frac{1}{[\boldsymbol{\Sigma}_{1}]_{z, z}
 +[\boldsymbol{\Sigma}_{2}]_{z,z}} \right)_{z \in \{1,\ldots,q\}} \right\} (\boldsymbol{\mu}_1-\boldsymbol{\mu}_2)\right].
\end{multline*}}

To get an $\varepsilon$-bracket for the   probability density functions, we have to construct a $\delta$-net for the variance and the mean \CG{for some value of} $\delta$ to be specified later.

\bigskip
\noindent 
\textit{Step 1: Construction of a net for the variance.}
Fix $\varepsilon \in (0,1]$ and {$\delta = \varepsilon/\sqrt{2}q$}.
Let $b_j = (1+\delta)^{1-{j/2}} A_{\boldsymbol{\Sigma}}$.
For each $j \in \{ 2, \ldots N\}$, we have $[a_{\boldsymbol{\Sigma}}, A_{\boldsymbol{\Sigma}}]= [b_N,b_{N-1}] \bigcup \cdots \bigcup [b_3,b_2]$, where $N$ is chosen to recover everything. We want
\begin{eqnarray*}
a_{\boldsymbol{\Sigma} }= (1+\delta)^{1-N/2} A_{\boldsymbol{\Sigma}}
&\Leftrightarrow&  \ln \frac{a_{\boldsymbol{\Sigma}}}{A_{\boldsymbol{\Sigma}}} = \left( 1-\frac{N}{2} \right) \ln(1+\delta)\\
&\Leftrightarrow& N = \frac{2 \ln( \frac{A_{\boldsymbol{\Sigma}}}{a_{\boldsymbol{\Sigma}}}\sqrt{1+\delta} )}{\ln (1+\delta)}.
\end{eqnarray*}
\CG{As we} want $N$ to be an integer, we choose 
$$
N= \left\lceil \frac{2 \ln( \frac{A_{\boldsymbol{\Sigma}}}{a_{\boldsymbol{\Sigma}}}\sqrt{1+\delta} )}{\ln (1+\delta)}  \right\rceil,
$$
\CG{where} $\lceil \cdot \rceil$ \CG{denotes} the ceiling function.
We \CG{then} get a regular net for the variance.
We could let $\mathbf B=\diag(b_{i(1)},\ldots, b_{i(q)})$, close to $\boldsymbol{\Sigma}$ (and deterministic, independent of the values of $\boldsymbol{\Sigma}$), where $i$ is a permutation 
such that $b_{i(z)+1} \leq [\boldsymbol{\Sigma}]_{z,z} \leq b_{i(z)}$ for all $z \in \{1,\ldots,q\}$.

\bigskip
\noindent 
\textit{Step 2: Construction of a net for the mean vectors}.
We use the singular decomposition of $\boldsymbol{\beta}$, $\boldsymbol{\beta} = \sum_{r=1}^{R} \sigma_r [\mathbf{u}]_{r,\bullet}{^\top} [\mathbf{v}]_{r,\bullet}$, with $(\sigma_r)_{1 \leq r \leq R}$ the singular values, 
and $([\mathbf{u}]_{r,\bullet})_{1 \leq r \leq R}$ and $([\mathbf{v}]_{r,\bullet})_{1\leq r \leq R}$ unit vectors.
\CG{These} vectors are also bounded.

We are looking for $\ell$ and $u$ such that $d_H( \ell,u) \leq \varepsilon$, and $\ell \leq f  \leq u$.
We  use a dilatation of a Gaussian to construct such an $\varepsilon$-bracket of $\varphi$.
We let
\begin{align*} 
\ell (\mathbf x,\mathbf y)&= (1+\delta)^{-(p^2qR+3q/4)} \varphi\{\mathbf y|\boldsymbol{\nu}_{J,R} \mathbf x,(1+\delta)^{-1/4} \mathbf{B}^{1}\} , \\
u (\mathbf x,\mathbf y)&= (1+\delta)^{p^2qR+3q/4} \varphi\{\mathbf y|\boldsymbol{\nu}_{J,R} \mathbf x, (1+\delta)\mathbf{B}^2\},
\end{align*}
where $\mathbf B^1$ and $\mathbf B^{2}$ are such that, for all $z \in \{1,\ldots,q\}$, $[\mathbf B^1]_{z,z} \leq \boldsymbol{\Sigma}_{z,z} \leq [\mathbf B^2]_{z,z}$ (see step $1$).

The means $\boldsymbol{\nu}_{J,R} \in \mathbb{R}^{q \times p}$  will be specified later. 
Just remark that $J$ is the set of the relevant columns, and $R$ the rank of $\boldsymbol{\nu}_{J,R}$: 
we  decompose $\boldsymbol{\nu}_{J,R} =\sum_{r=1}^\mathbf{R} \tilde{\sigma}_r [\mathbf{\tilde{u}}]_{r,\bullet}{^\top} [\mathbf{\tilde{v}}]_{r,\bullet}$,
with $\mathbf{ \tilde{u}} \in \mathbb{R}^{R \times |J|}$ and $\mathbf{\tilde{v}} \in \mathbb{R}^{R \times q}$. We have 
$$
\ell (\mathbf x,\mathbf y) \leq f(\mathbf y|\mathbf x) \leq u(\mathbf x,\mathbf y) \quad
\Leftrightarrow \quad
||\boldsymbol{\beta} \mathbf x -\boldsymbol{\nu}_{J,R} \mathbf x ||_2^2 \leq p^2qR \frac{\delta^2}{2} a_{\boldsymbol{\Sigma}}^2(1-2^{-1/4}).
$$
Note that $||\boldsymbol{\beta} \mathbf x -\boldsymbol{\nu}_{J,R} \mathbf x||_2^2 \leq p ||\boldsymbol{\beta}-\boldsymbol{\nu}_{J,R}||_2^2 ||\mathbf x||_\infty$. \CG{Thus we need}
\begin{align}
\label{controleNorme2}
||\boldsymbol{\beta}-\boldsymbol{\nu}_{J,R}||_2^2 \leq pqR \, \frac{\delta^2}{2} \, a_{\boldsymbol{\Sigma}}^2 (1-2^{-1/4}). 
\end{align}

According to \cite{inegOracleLassoMLE}, $d_H(\ell,u) \leq 2(p^2qR +3q/4)^2 \delta^2$ \CG{and hence by choosing
$\delta = {\varepsilon}/\{\sqrt{2}(pqR+3/4q)\}$, we get the desired bound. Let us now see how} to construct $\boldsymbol{\nu}_{J,R}$ to get \eqref{controleNorme2}. First write
\begin{eqnarray*}
||\boldsymbol{\beta}-\boldsymbol{\nu}_{J,R}||_2^2 &=& \sum_{j=1}^p \sum_{z=1}^q \left| 
 \sum_{r=1}^{R} \sigma_r [\mathbf{ {u}}]_{r,j} [\mathbf{ v}]_{r,z} - \tilde{\sigma}_r [\mathbf{ \tilde{u}}]_{r,j} [\mathbf{\tilde{ v}}]_{r,z}\right|^2\\
&=& \sum_{j,z} |\sum_{r=1}^{R}  
 |\sigma_r-\tilde{\sigma}_r||[\mathbf{u}]_{r,j}[\mathbf{ v}]_{r,z}| \\
 && - \tilde{\sigma}_r |[\mathbf{ \tilde{u}}]_{r,j}-[\mathbf{ {u}}]_{r,j}||[\mathbf{\tilde{ v}}]_{r,z}| 
 - \tilde{\sigma}_r [\mathbf{ {u}}]_{r,j} |[\mathbf{ v}]_{r,z}-[\mathbf{\tilde{ v}}]_{r,z}| |^2.
\end{eqnarray*}
Then
\begin{eqnarray*}
||\boldsymbol{\beta}-\boldsymbol{\nu}_{J,R}||_2^2 &\leq& \sum_{j=1}^p \sum_{z=1}^q  \left| \sum_{r=1}^{R}|\sigma_r-\tilde{\sigma}_r|+A_\sigma |[\mathbf{ \tilde{u}}]_{r,j}-[\mathbf{ {u}}]_{r,j}| + A_\sigma |[\mathbf{ v}]_{r,z}-[\mathbf{\tilde{ v}}]_{r,z}| \right|^2\\
 &\leq& 2pqR (\max_r |\sigma_r-\tilde{\sigma}_r|^2 + A_\sigma \max_{r,j} |[\mathbf{ \tilde{u}}]_{r,j}-[\mathbf{ {u}}]_{r,j}|^2 \\
&& \hspace{1.6cm}+A_\sigma \max_{r,z} |[\mathbf{\tilde{ v}}]_{r,z}-[\mathbf{ v}]_{r,z}|^2 ).
\end{eqnarray*}
If we now choose $\tilde{\sigma}_r$, $[\mathbf{ \tilde{u}}]_{r,j}$ and $[\mathbf{\tilde{ v}}]_{r,z}$ such that 
\begin{align*}
|\sigma_r-\tilde{\sigma}_r| &\leq \frac{\delta}{\sqrt{12}} \, a_{\boldsymbol{\Sigma} }\sqrt{1-2^{-1/4}},\\
|[\mathbf{ {u}}]_{r,j}-[\mathbf{ \tilde{u}}]_{r,j}| &\leq \frac{\delta}{\sqrt{12}A_\sigma} \, a_{\boldsymbol{\Sigma} }\sqrt{1-2^{-1/4}},\\
|[\mathbf{ v}]_{r,z}-[\mathbf{\tilde{ v}}]_{r,z}| &\leq \frac{\delta}{\sqrt{12}A_\sigma} \, a_{\boldsymbol{\Sigma} }\sqrt{1-2^{-1/4}},
\end{align*}
\CG{we can conclude, as desired, that}
$$
||\boldsymbol{\beta}-\boldsymbol{\nu}_{J,R}||_2^2 \leq pqR \, \frac{\delta^2}{2} \, a_{\boldsymbol{\Sigma}}^2 (1-2^{-1/4}).
$$
\CG{To accomplish this, let $\lfloor \cdot \rfloor$ the floor function and, for all $r \in \{1,\ldots,R\}, j \in \{1,\ldots,p\}, z \in \{1,\ldots,q\}$, set}
\begin{align*}
S &= \mathbb{Z} \cap \left[ 0, \left\lfloor \frac{A_\sigma}{\frac{\delta}{\sqrt{12}} a_{\boldsymbol{\Sigma} }\sqrt{1-2^{-1/4}}}\right\rfloor\right], \\
\tilde{\sigma}_r &=  \underset{\varsigma \in S}{\operatorname{argmin}} \left|\sigma_r - \frac{\delta}{\sqrt{12}} a_{\boldsymbol{\Sigma} }\sqrt{1-2^{-1/4}} \varsigma \right| , \\
 U &= \mathbb{Z} \cap \left[ 0, \left\lfloor \frac{A_\sigma}{\frac{\delta}{\sqrt{12}A_\sigma} a_{\boldsymbol{\Sigma} }\sqrt{1-2^{-1/4}}}\right\rfloor\right] , \\
[\mathbf{ \tilde{u}}]_{r,j} &=  \underset{\mu \in U}{\operatorname{argmin}} \left|[\mathbf{ {u}}]_{r,j} - \frac{\delta}{\sqrt{12}A_\sigma} a_{\boldsymbol{\Sigma} }\sqrt{1-2^{-1/4}} \mu \right| , \\
 V& = \mathbb{Z} \cap \left[ 0, \left\lfloor \frac{A_\sigma}{\frac{\delta}{\sqrt{12}A_\sigma} a_{\boldsymbol{\Sigma} }\sqrt{1-2^{-1/4}}}\right\rfloor\right] , \\
[\mathbf{\tilde{ v}}]_{r,z} &=  \underset{\nu \in V}{\operatorname{argmin}} \left|[\mathbf{ v}]_{r,z} - \frac{\delta}{\sqrt{12}A_\sigma} a_{\boldsymbol{\Sigma} }\sqrt{1-2^{-1/4}} \nu \right| .
\end{align*}

\CG{Note} that we just need to determine vectors $(([\mathbf{ \tilde{u}}]_{r,j})_{1\leq j \leq J-r})_{1 \leq r \leq R}$ and 
$(([\mathbf{\tilde{ v}}]_{r,z})_{1\leq z \leq q-r})_{1 \leq r \leq R}$ because \CG{these are unit vectors}.
Then we let, for all $z \in \{1,\ldots,p\}$,
$$
[\boldsymbol{\nu}_{J,R}]_{z,j} = \begin{cases} 
\sum_{r=1}^{R} \tilde{\sigma}_r [\mathbf{ \tilde{u}}]_{r,j} [\mathbf{\tilde{ v}}]_{r,z} & \mbox{if } j \in J, \\
0 & \mbox{if } j \in J^c.
\end{cases}
$$

\bigskip
\noindent
\textit{Step 3: Upper bound of the number of $\varepsilon$-brackets for $\mathcal{F}_{(J,{R})}$.}
We have defined the brackets.
Let $c=(1-2^{-1/4})/12$.
We want to control the entropy.
We denote by $\mathcal{B}_\varepsilon (\mathcal{F}_{(J,{R})} )$ the minimum number of brackets of width $d_H^{\otimes_n}$ smaller than $\varepsilon$ such that every density of $\mathcal{F}_{(J,{R})}$ belongs to one of those brackets. We have
\begin{align*}
 |\mathcal{B}_\varepsilon (\mathcal{F}_{(J,{R})} )| & \leq \sum_{l=2}^N \left( \frac{A_\sigma}{\delta a_{\boldsymbol{\Sigma} }\sqrt{c}} \right)^{R} \left( \frac{A_\sigma^2}{\delta a_{\boldsymbol{\Sigma} }\sqrt{c}} \right)^{R\left( \frac{2J-R-1}{2}+\frac{2q-R-1}{2}\right)}\\
&  \leq(N-1) \left( \frac{A_\sigma^2}{\delta a_{\boldsymbol{\Sigma} }\sqrt{c}} \right)^{R(J+q-R)} A_\sigma^{-R} \\
&  \leq  C(a_{\boldsymbol{\Sigma}},A_{\boldsymbol{\Sigma}},A_\sigma,J,{R}) \delta^{-D_{(J,{R})}-1}
 \end{align*}
with 
$$
C(a_{\boldsymbol{\Sigma}},A_{\boldsymbol{\Sigma}},A_\sigma,J,{R}) = 2 \left( \frac{A_{\boldsymbol{\Sigma}}}{a_{\boldsymbol{\Sigma}}}+\frac{1}{2}\right)\left( \frac{A_\sigma^2}{ a_{\boldsymbol{\Sigma} }\sqrt{c}} \right)^{R(J+q-R)} A_\sigma^{-R}.
$$

\bigskip
\noindent
\textbf{B.3 For the mixture}. 
We begin by computing the bracketing entropy: according to \eqref{entropie},
\begin{align*}
 \mathcal{H}_{[\cdot]} (\varepsilon, \mathcal{S}^{\mathcal{B}}_{(K,J,\mathbf{R})},d_H) &\leq \ln \left\{C \left( \frac{1}{\varepsilon}\right)^{D_{(K,J,\mathbf{R})}}\right\},
\end{align*}
where
\begin{multline}
\label{constante}
C = 2^K K(2\pi e)^{K/2} 3^{K-1} \left( \frac{A_{\boldsymbol{\Sigma}}}{a_{\boldsymbol{\Sigma}}} +\frac{1}{2} \right)^K 
\\ \times \left( \frac{A_\sigma \sqrt{12}}{a_{\boldsymbol{\Sigma}}^2 \sqrt{1-2^{-1/4}} } \right)^{D_{(K,J,\mathbf{R})}} A_\sigma^{-\sum_{k=1}^K \mathbf{R}_k},
\end{multline}
and $D_{(K,J,\mathbf{R})}= \sum_{k=1}^K \mathbf{R}_k (|J|+q-\mathbf{R}_k)$. We \CG{must} determine $\psi_{(K,J,\mathbf{R})}$ such that 
\begin{align}
\label{2etoiles}
\int _0^\varpi \sqrt{ \mathcal{H}_{[\cdot]} (\varepsilon,\mathcal{S}^{\mathcal{B}}_{(K,J,\mathbf{R})}(s^{(K,J,\mathbf{R})},\varpi),d_H^{\otimes_n})} d \varepsilon \leq \psi_{(K,J,\mathbf{R})}(\varpi). 
\end{align}
\CG{To this end, we} compute the integral, viz.
\begin{align*}
& \hspace{-2cm} \int _0^\varpi \sqrt{ \mathcal{H}_{[\cdot]} (\varepsilon,\mathcal{S}^{\mathcal{B}}_{(K,J,\mathbf{R})}(s^{(K,J,\mathbf{R})},\varpi),d_H^{\otimes_n})} d \varepsilon \\
&\leq \varpi \sqrt{\ln(C)} + \sqrt{D_{(K,J,\mathbf{R})}} \int_{0}^\varpi \sqrt{\ln\left( \frac{1}{\varepsilon} \right)} d\varepsilon\\
&\leq \sqrt{D_{(K,J,\mathbf{R})}} \varpi \left\{ \sqrt{\pi}+\sqrt{\frac{\ln(C)}{D_{(K,J,\mathbf{R})}}} +\sqrt{\ln{\left( \frac{1}{\varpi \wedge 1} \right) }} \right\}
 \end{align*}
with, according to \eqref{constante},
\begin{align*}
  \ln(C) & =  K \ln(2) +\frac{K}{2} \ln(2\pi e) +(K-1)\ln(3) +K \ln\left(\frac{A_{\boldsymbol{\Sigma}}}{a_{\boldsymbol{\Sigma}}}+\frac{1}{2}\right) \\
  &+D_{(K,J,\mathbf{R})} \ln\left( \frac{A_\sigma^2 \sqrt{12}}{a_{\boldsymbol{\Sigma} }\sqrt{1-2^{-1/4}}}\right) +\ln(K) +\sum_{k=1}^K \mathbf{R}_k \ln \left(\frac{1}{A_\sigma} \right)\\
& \leq D_{(K,J,\mathbf{R})} \left\{ \ln(2) +\ln(2 \pi e)+ \ln{3} +1 +\ln\left( \frac{A_\sigma^2 \sqrt{12}}{a_{\boldsymbol{\Sigma} }\sqrt{1-2^{-1/4}}}\right)\right\}\\
&+D_{(K,J,\mathbf{R})} \left\{\ln\left( \frac{A_{\boldsymbol{\Sigma}}}{a_{\boldsymbol{\Sigma}}}+\frac{1}{2} \right)  \right\} \\
&   \leq D_{(K,J,\mathbf{R})} \left[ \ln\left( \frac{12\sqrt{12} \pi e}{\sqrt{1-2^{-1/4}}}\right) + \ln\left\{\frac{A^2_\sigma}{a_{\boldsymbol{\Sigma}}} 
  \left(\frac{A_{\boldsymbol{\Sigma}}}{a_{\boldsymbol{\Sigma}}}+\frac{1}{2} \right)\right\} \right].
\end{align*}
Then
\begin{multline*}
 \int _0^\varpi \sqrt{ \mathcal{H}_{[\cdot]} (\varepsilon,\mathcal{S}^{\mathcal{B}}_{(K,J,\mathbf{R})}(s^{(K,J,\mathbf{R})},\varpi),d_H^{\otimes_n})} d \varepsilon \\
 \leq \sqrt{D_{(K,J,\mathbf{R})}} \left\{ 2+ \sqrt{\ln \left\{ \frac{A_\sigma^2}{a_{\boldsymbol{\Sigma}}}\left( \frac{A_{\boldsymbol{\Sigma}}}{a_{\boldsymbol{\Sigma}}}+\frac{1}{2} \right) 
 \right\}} + \sqrt{\ln\left(\frac{1}{\varpi \wedge 1 } \right) }\right\}.
\end{multline*}
Consequently, by putting 
$$
B=2+ \sqrt{\ln \left\{ \frac{A_\sigma^2}{a_{\boldsymbol{\Sigma}}}\left( \frac{A_{\boldsymbol{\Sigma}}}{a_{\boldsymbol{\Sigma}}}+\frac{1}{2} \right)\right\}} ,
$$
we find that the function $\psi_{(K,J,\mathbf{R})}$ defined on \CG{$(0, \infty)$} by
$$
\psi_{(K,J,\mathbf{R})} (\varpi) = \sqrt{D_{(K,J,\mathbf{R})}} \varpi \left\{ B+\sqrt{ \ln \left(\frac{1}{\varpi \wedge 1} \right)} \right\}
$$
satisfies \eqref{2etoiles}.
Besides, $\psi_{(K,J,\mathbf{R})}$ is non-decreasing and $\varpi \mapsto \psi_{(K,J,\mathbf{R})}(\varpi) /\varpi$ is non-increasing, \CG{so} $\psi_{(K,J,\mathbf{R})}$ is convenient.

Finally, we need to find an upper bound \CG{on} $\varpi_*$ satisfying 
$\psi_{(K,J,\mathbf{R})} (\varpi_*) = \sqrt{n} \varpi^2_*$. Consider $\varpi_*$ such that $\psi_{(K,J,\mathbf{R})} (\varpi_*) = \sqrt{\varpi_*^2}$. This is equivalent to \CG{solving}
$$
\varpi_*= \sqrt{\frac{D_{(K,J,\mathbf{R})}}{n}} \left\{ B+ \sqrt{\ln\left(\frac{1}{\varpi_* \wedge 1} \right) } \right\}.
$$
We then choose
$$\varpi_*^2 \leq \frac{D_{(K,J,\mathbf{R})}}{n}\left\{ 2B^2 + \ln\left(\frac{1}{1 \wedge \frac{D_{(K,J,\mathbf{R})}}{n} B^2 } \right) \right\}.$$

\bigskip
\noindent
\textbf{C Assumption \ref{AsumptionK}}. 
Let 
$$w_{(K,J,\mathbf{R})}=\tilde{D}_{(K,J)} \ln\left\{ \frac{4epq}{(\tilde{D}_{(K,J)}-q^2)\wedge pq} \right\} +\sum_{k \in \{1,\ldots,K\}}\mathbf{R}_k,$$ 
where $\tilde{D}_{(K,J)} = K(1+|J|)$.
Then
\begin{align*}
 \sum_{(K,J,\mathbf{R})}e^{-w_{(K,J,\mathbf{R})}}& \leq \sum_{K \geq 1} 
 \left( \sum_{R \geq 1} e^{-R} \right)^K  \left\{ \sum_{1 \leq |J| \leq p} e^{-\tilde{D}_{(K,J)} \ln\left( \frac{4epq}{(\tilde{D}_{(K,J)}-q^2)\wedge pq} \right)} \right\} \\
 & \leq \sum_{K \geq 1} 1^K \left\{ \sum_{1 \leq |J| \leq p} e^{-\tilde{D}_{(K,J)} \ln\left( \frac{4epq}{(\tilde{D}_{(K,J)}-q^2)\wedge pq} \right)} \right\} \le 2,
\end{align*}
\CG{where the last inequality stems from Proposition~4.5 in \cite{inegOracleLassoMLE}}.

We \CG{observe} that the weights depend on $(K,J,\mathbf{R})$, which leads to a \CG{sharper} control \CG{but forces} the penalty to depend on $\tau$, the parameter related to \CG{Bernstein's} inequality, according to \eqref{penalite2}. \CG{Larger weights can be used to avoid this dependence, however}. We prefer to use \CG{the present weights} to get the thinest control possible.

\end{document}